\documentclass[12pt,a4paper,twoside,reqno]{amsart}

\usepackage[latin1]{inputenc}

\usepackage{graphicx}
\usepackage{subfigure}

\usepackage{amsmath}
\usepackage{amsfonts}
\usepackage{amssymb}
\usepackage{amstext}
\usepackage{amsbsy}
\usepackage[active]{srcltx}

\newtheorem{corollary}{\sc Corollary}[section]
\newtheorem{theorem}{\sc Theorem}[section]
\newtheorem{lemma}{\sc Lemma}[section]
\newtheorem{proposition}{\sc Proposition}[section]
\newtheorem{remark}{\sc Remark}
\newtheorem{definition}{\sc Definition}%[section]

\newtheorem{example}{\sc Example}[section]

\newcommand{\thmref}[1]{Theorem~\ref{#1}}

\newcommand{\lemref}[1]{Lemma~\ref{#1}}

\newcommand{\coref}[1]{Corollary~\ref{#1}}

\newcommand{\defref}[1]{Definition~\ref{#1}}
\newcommand{\exaref}[1]{Example~\ref{#1}}

\let\hip=\Hip

\def \Cx {{\mathbb C}}
\def\Na{\mathbb N}

\newcommand{\R}{\mathbb R}

\def\De{\Delta}

\def \Rea {{\mathbb R}}

\def\Gai{\Ga_infty}

\def\hiper{\mathbb H}
\let\h=\hip
\let\hipe=\hiper

\let\re=\R

\def\Mc{\mathcal M}

\def\Sc{\mathcal S}

\def\Nc{\mathcal N}
\def\Dc{\mathcal D}

\def\s{\mathbb  S}

\def\Om{\Omega}

\def \om{\omega}
\def \ga{\gamma}

\def \be{\beta}

\def \ga{\gamma}
\def \Ga{\Gamma}

\def \e{\varepsilon}

\def\Gai{\Ga_\infty}

\def \ep{\epsilon}

\def\vphi{\varphi}

\def \om{\omega}
\def \Om{\Omega}

\def\leqs{\leqslant}
\def\geqs{\geqslant}

\def\rmd{\mathop{\rm d\kern -1pt}\nolimits}
\def\rme{\mathop{\rm e\kern -1pt}\nolimits}

\def\cro{\Mc_\rho}

\def\e{\varepsilon}

\DeclareMathOperator{\diver}{div}
\DeclareMathOperator{\dist}{dist}

\DeclareMathOperator{\ext}{ext}

\DeclareMathOperator{\inter}{int}

\def \Cc{{\mathcal C}}

\def\bel{ \medskip
 \centerline{$ \ast \hbox to 1.0cm{}\ast \hbox to 1.0cm{}\ast $}
}

\def\overl{\overline}

\def\goto{\rightarrow}

\def\do{\downarrow}

\def\longerrightarrow{-\kern-5pt\longrightarrow}

\def\vphi{\varphi}
\def\star{\lower 1pt\hbox{*}}
\def \nulset {
\raise 1pt\hbox{ \hskip -3pt$\not$\kern -0.2pt \raise
.7pt\hbox{${\scriptstyle\bigcirc}$}}}

\def \ep{\epsilon}

\newcommand{\hd}{\mathbb{H}^2}

\newcommand{\hi}[1]{\mathbb{H}^#1}
\newcommand{\m}[1]{\mathbb{R}^#1}
\newcommand{\ch}{\cosh}

\newcommand{\pain}{\partial_{\infty}}

\newcommand{\ov}[1]{\overline{#1}}

\def\PO{\partial\Om\cup\partial_\infty\Om}

\let\leq=\leqslant
\let\geq=\geqslant

\let\dc=\Dc

\begin{document}

\title[minimal]
{An asymptotic theorem for minimal surfaces 
and existence results for
minimal graphs in $\hip^2 \times \R$}

\author[  R. Sa Earp  and E. Toubiana]
{\scshape R. Sa Earp and E. Toubiana}

 \address{Departamento de Matem\'atica,
  Pontif\'\i cia Universidade Cat\'olica do Rio de Janeiro, Rio de Janeiro,  22453-900 RJ,
 Brazil }\email{earp@mat.puc-rio.br}

\address{Institut de Math{\'e}matiques de Jussieu, Universit{\'e} Paris
VII, Denis Diderot, Case 7012,
         2 Place Jussieu,
         75251 Paris Cedex 05, France}
\email{toubiana@math.jussieu.fr}
\thanks{The authors would like to thank CNPq,  PRONEX of
Brazil and Accord Brasil-France, for partial financial support}

\date{\today}


\maketitle

\begin{abstract}

In this paper 
we prove a general and sharp Asymptotic Theorem
for minimal surfaces in $\hd \times \R$.
As a consequence, we prove that there is no properly immersed  minimal
surface whose asymptotic boundary $\Gamma_\infty$ is a Jordan curve
homologous to zero in $\partial_\infty\hip^2\times \R$ such that
$\Gamma_\infty$ is contained in a  slab between two horizontal circles
of $\partial_\infty\hip^2\times \R$  with width  equal to
$\pi.$

We construct minimal  vertical graphs in
$\hip^2\times \R$ over certain unbounded admissible domains taking
certain prescribed finite boundary data and certain prescribed
asymptotic boundary data. Our admissible unbounded domains $\Om$
in $\hip^2\times \{0\}$  are non necessarily convex and non
necessarily bounded by convex arcs; each component of its boundary
is properly embedded with zero, one or two points on its
asymptotic boundary, satisfying a further geometric condition.
\end{abstract}

\section{Introduction }

In this paper we prove an Asymptotic Theorem for minimal 
surfaces in $\hd \times \R$.
 Indeed, we prove a
 surprising general and sharp nonexistence result. As a consequence, we deduce 
that
there is no complete  properly immersed minimal surface whose asymptotic
boundary
$\Ga_\infty$ is a Jordan curve homologous to zero in
$\partial_\infty\hip^2\times \R$ contained in an
open slab between two horizontal circles of
$\partial_\infty\hip^2\times \R$  with width equal 
to $\pi.$ The last statement is still true in a closed slab with width 
{\em equal} to  $\pi$ in the class of minimal surfaces continuous up to the
asymptotic boundary.
This result is sharp  in the following sense. 
We show that for any $\ell >\pi$ there is a Jordan curve  
$\Ga_\infty \subset \pain\hd \times \R$ homologous to zero with vertical height 
equal to $\ell$ which is the asymptotic boundary of a complete minimal surface,
continuous up to its asymptotic boundary $\Ga_\infty$. Moreover, this surface
is 
 invariant by
hyperbolic translations and is constituted of two minimal vertical graphs
over the exterior of an equidistant curve, symmetric about the
horizontal slice $\hip^2\times \{0\}$.
In fact, 
the Jordan curve  $\Ga_\infty$ is the union of two vertical segments with
two half circles in $\partial_\infty \hip^2$. Another consequence of our
Asymptotic Theorem is that there is no complete properly immersed minimal
surface contained in an open slab of width equal to $\pi$ of $\hd \times \R$,
such that the vertical projection of its asymptotic boundary on $\pain \hd
\times \{0\}$ omits an open arc.

 Those results contrast with an
analogous situation when the ambient space is the hyperbolic
three-space $\hip^3,$ due to the existence of the minimal vertical
graph (taking the upper half-space model) whose asymptotic
boundary is any convex curve lying in $\pain \hip^3.$ Indeed, the
authors have solved the Dirichlet Problem in
$\hip^3$ for the minimal vertical equation over a convex domain
$\Om$ in  $\pain \hip^3,$ taking any prescribed continuous
boundary data on $\partial \Om$ (\cite{asi}). There are also the
general results proved by M. Anderson ~\cite{An1} and ~\cite{An2}.

We give some geometric conditions  to
 construct minimal  vertical graphs in
$\hip^2\times \R$ over certain unbounded admissible domains taking
certain prescribed finite boundary data and certain prescribed
asymptotic boundary data.

To obtain our existence results we establish the {\em Perron
process} when the finite boundary data and the asymptotic boundary
data are continuous except maybe at a finite set. 

As a consequence, we prove the following. Let
$\Om$ be a convex unbounded
domain. Let $g:\partial \Omega \cup \pain \Omega\rightarrow \R$ be a bounded
function everywhere continuous except at a finite set $S$. Then  
$g$ admits an extension $u$ satisfying the minimal vertical
equation over $\Om$ such that the total boundary of the graph of $u$ is the
union of the graph of $g$ on 
$(\partial \Omega \cup \pain \Omega)\setminus S$
 with vertical segments at the points of $S$. This result was obtained
independently by M. Rodr\' \i guez and H. Rosenberg

We built {\em barriers
 } at each convex point of a convex finite boundary, where the boundary data are
continuous and bounded and we construct  {\em barriers} 
at each point of the asymptotic boundary 
where the
asymptotic data is continuous.
 Our admissible unbounded domains $\Om$ in $\hip^2\times \{0\}$  are non necessarily convex and non
necessarily bounded by convex arcs; each component of its boundary
is properly embedded with zero, one or two points on its
asymptotic boundary, satisfying a further geometric condition:
each connected component $C_0$ of $\partial \Om$ satisfies the
{\em Exterior circle of $($uniform$)$ radius $\rho$ condition}.
Particularly, we consider an admissible domain $\Om$ that is the
exterior of a $C^2$ Jordan curve $\Ga$ in the horizontal slice.

 We obtain the existence of minimal graph $M$ over an admissible
domain $\Om$
 in $\hip^2\times \{0\}$ such
that  the finite boundary of $M$ is $\partial\Om$ and the asymptotic
boundary of $M$ is a certain Jordan curve $ \Ga_\infty$ constituting
in the union of
 bounded continuous vertical graphs with  the vertical segments joining the
points of
discontinuities, such that $ \Ga_\infty$ is contained inside a
certain slab of $\pain\hip^2\times \R$ depending on the geometry of $\Om.$

We consider admissible domains, that we call E-admissible domains,
such that each component of the boundary has two points at
its asymptotic boundary and has at each point of its finite
boundary an exterior equidistant curve.  We obtain  analogous
existence results for E-admissible domains.

\section{An asymptotic theorem}

In this section we prove an Asymptotic Theorem that ensures 
some nonexistence results about minimal surfaces with
some given asymptotic boundary. 

\begin{theorem}[Asymptotic Theorem]
\label{T.asymp.boundary}$ $

 Let $\ga\subset \pain \hd \times \R$ be an arc. Assume there exist a vertical
straight line $L\subset \pain \hd \times \R$ and a subarc 
$\ga^\prime \subset \ga$ such that
\begin{enumerate}
\item  $\ga^\prime \cap L\not= \emptyset$ and 
$\partial \ga^\prime \cap L=\emptyset$,
\item  $\ga^\prime$ stays on one side of $L$,
\item $\ga^\prime \subset \pain \hd \times \mathopen (t_0,\pi +t_0)$, for some
real number $t_0$.
\end{enumerate}
Therefore, there is no properly immersed minimal surface 
$($maybe with finite boundary$)$,
$M\subset \hd \times\R$, with asymptotic boundary
$\ga$ and such that $M\cup \ga$ is a continuous surface
with boundary.
\end{theorem}

\begin{proof}
By assumption there exists a point $p$ in $\ga^\prime \cap L$. If there is a
vertical segment in $\ga^\prime \cap L$, we choose $p$ to be the midpoint of
this segment. Up to a vertical
translation, we can assume that $p\in \pain \hd \times \{0\}$. The vertical
projection of $\ga^\prime$ on $\pain \hd \times \{0\}$ is an arc $\beta$ with
$p$ as one of the two end points. Let $\e >0$ be a real number
to be chosen later. Let $q_1,q_2\in \pain\hd \times \{0\}$ be two distinct
points  such that $q_1\in \beta$, $q_2\not\in \beta$ and the Euclidean distance
on $\pain\hd \times \{0\}$ from $p$ to $q_i$ is $\e$, $i=1,2$. Let 
$c\subset \pain \hd \times \{0\}$ be the complete geodesic with asymptotic
boundary $\{q_1,q_2\}$ and let $S=c\times \R$ be the vertical geodesic plane
defined by $c$. 

 Let $M\subset \hd \times \R$ be a minimal surface (if any) with asymptotic
boundary $\ga$ and such that $M\cup \ga$ is a continuous surface with boundary.
If $\e$ is small enough we have $S\cap \partial M=\emptyset$.
Let $M_0\subset M$ be the connected component of $M\backslash S$ 
containing $p $ in
its asymptotic boundary. Therefore, the asymptotic boundary of $M_0$ is
a subarc $\ga_0$ of $\ga ^\prime$ containing $p$ in its interior: 
$\ga_0 \subset \ga^\prime \subset \ga$ and $p\in \textrm{Int}(\ga_0)$. 
Let $\beta^\prime \subset \beta \subset \pain \hd \times \{0\}$ be the subarc of
$\beta$ with end points  $p$ and $q_1$. For $\e$ small
enough we have $\ga_0 \subset \beta^\prime \times \mathopen (-\pi/2,
\pi/2\mathclose)$. By construction there exist two real numbers $a$ and $b$
satisfying $a<0<b$, $b-a<\pi$ and $\partial \ga_0=\{(q_1,a), (q_1,b)\}$. 

 Observe that, by continuity, for $\e$ small enough the whole component $M_0$ is
inside the slab $\hd \times \mathopen (-\pi/2, \pi/2\mathclose)$.
Furthermore, the finite boundary $\partial M_0$ of $M_0$ is contained in the
vertical geodesic plane $S$. Therefore, there is a complete geodesic 
$c_1\subset \hd \times \{0\}$ with asymptotic boundary in the open arc 
$\mathopen (p,q_2 \mathclose )\subset \pain\hd \times \{0\} \backslash \beta$,
such that $M_0 \cap (c_1\times \R)=\emptyset$.

  Let $C\subset \hd\times \R$ be a complete catenoid  whose a component of the
asymptotic boundary stays at height $T_1$ and the other component at height
$T_2$ such that $T_1 <a <b <T_2$ (such a catenoid exists since $0<b-a<\pi$),
note that $T_2-T_1<\pi$, the reader can see the geometric behaviour of the
catenoids in Lemma \ref{l1} or \cite{N-SE-S-T}. By continuity, we can choose
$T_1$ and $T_2$ such that $M_0$ is entirely contained in the open slab
$\hd \times  (T_1, T_2)$.
Finally, let $t_1,t_2$ be two real numbers satisfying \newline
$T_1<t_1<a <b<t_2<T_2$, such that $M_0$ is entirely contained in the open slab
$\hd \times  (t_1, t_2)$. 

Let $\ov C$ be the part of $C$ contained in the
slab
$\{t_1 \leqs t\leqs t_2\}$, that is, $\ov C=C\cap (\hd \times [t_1,t_2])$.
Observe that $\ov C$ is a compact surface. Up to a hyperbolic translation we
can send $\ov C$ into the connected component of 
$\hd \times \R \backslash (c_1\times \R)$ not containing $p$ in its asymptotic
boundary, so we can assume that $\ov C$ has this property. 

 Let $c_2\subset \hd \times \{0\}$ be a complete geodesic with an asymptotic
boundary point in $\textrm{Int} (\beta^\prime)$ and the other asymptotic
boundary point 
in the open arc $\mathopen ( p,q_2\mathclose ) $. We choose $c_2$ such that
$c_1$ is contained in the component of $(\hd \times \{0\})\setminus c_2$
containing $p$ in its asymptotic boundary.
 Consider the hyperbolic translations along $c_2$. Observe that all translated
copies of $\ov C$ have a component of the finite boundary at height $t_1$ and
the other component at height $t_2$. Therefore the boundary of any translated
copy of $\ov C$ has no intersection with $\ov{ M_0}$. Consequently some
translated
copy
of $\ov C$ must achieve a first interior contact point with $M_0$, which
contradicts
the maximum principle. This concludes the proof of the Theorem.

\end{proof}

\begin{corollary}\label{C.asym.bound}
Let $\Gamma_\infty \subset \pain \hd \times \R$ be a Jordan curve
homologous to zero $($in $\pain \hd \times \R)$. We have the following:
\begin{enumerate}
 \item \label{Item.less.pi} Suppose that  $\Gamma_\infty $ is strictly contained
in a 
closed slab between two horizontal circles of
$\partial_\infty\hip^2\times \R$  with width  equal
to $\pi$. Then,

\begin{enumerate}

\item \label{Item.Item.boundary} there is no properly immersed  minimal
surface $M$ with asymptotic boundary $\Gamma_\infty $, possibly with finite
boundary, 
such that $M\cup \Gamma_\infty$ is a continuous
surface with boundary.

\item \label{Item.Item.complete} there
is no complete properly immersed  minimal surface 
with asymptotic boundary
$\Gamma_\infty$ $($without any  assumption on $M\cup \Gamma_\infty )$.

\end{enumerate}

\item \label{Item.equal.pi}Suppose that $\Gamma_\infty $ is contained in a
slab with 
width  equal to $\pi$ but is not contained in any slab with
width  strictly less than $\pi$. Then, there is no 
complete minimal surface properly immersed in $\hd \times \R$, with asymptotic
boundary $\Gamma_\infty$ such that $M\cup \Gamma_\infty$ is a continuous surface
with boundary.

\end{enumerate}

\end{corollary}

\begin{proof}
 The Statement (\ref{Item.Item.boundary})
is a direct consequence
of Theorem \ref{T.asymp.boundary}. The Statement
(\ref{Item.Item.complete})
is a direct consequence of the proof
of Theorem \ref{T.asymp.boundary}.

Let us prove the Statement (\ref{Item.equal.pi}). Assume
there exists a properly immersed complete
minimal surface
$M$ with asymptotic boundary $\Gamma_\infty \subset \pain\hd \times [0,\pi]$.
By the maximum principle, we
deduce that $\textrm{Int}(M)  \subset \hd \times (0,\pi)$.
We deduce from Theorem \ref{T.asymp.boundary} that $\Gamma_\infty$ is
constitued by two vertical segments of length $\pi$: $\{q_1\}\times[0,\pi]$ and
$\{q_2\}\times[0,\pi]$, $q_i \in \pain\hd$ (identified with 
$\pain \hd \times \{0\}$), $i=1,2$, and two simple arcs 
$c,\ga \subset \pain\hd \times [0,\pi]$, the arc $c$ joining the points
$(q_1,0)$ and $(q_2,0)$ and the arc $\ga$ joining the points $(q_1,\pi)$
and $(q_2,\pi)$. Therefore we have
$$
\Gamma_\infty=(\{q_1\}\times[0,\pi])\cup\ga \cup (\{q_2\}\times[0,\pi])\cup c.
$$
Up to an ambiant isometry, we can assume that $q_1=e^{i\pi/4}$,
$q_2=e^{-i\pi/4}$ and that the vertical projection of $\Gamma_\infty$ on 
$\pain \hd \times \{0\}$ is the arc \newline
$\{(e^{i\theta},0)\mid -\pi/4\leqs \theta \leqs \pi/4\} $.

 Let $H$ be the parabolic complete minimal surface (foliated
by horocycles) whose asymptotic boundary is the vertical segment 
$\{-1\}\times [0,\pi]\subset \pain \hd \times \R$ with 
$(\pain \hd \times\{0\}) \cup (\pain \hd \times\{\pi\})$, see \cite{Daniel},
\cite{Haus} and \cite{Sa}. The ``neck'' of $H$  is a
horocycle $N$ in the slice $\hd \times \{\pi/2\}$.

\bigskip

{ \em Claim.} If $N\cap M=\emptyset$ then $H\cap M=\emptyset$.

\medskip

 Assume by contradiction that $N\cap M=\emptyset$ and 
$H^+\cap M\not=\emptyset$, where $H^+:=H\cap(\hd \times [\pi/2,\pi])$.  
For any $\e>0$ we denote by $H_\e^+$ the $\e$-vertical translated of $H^+$: 
$H_\e^+= H^+ +\e \partial/\partial t$.
Observe that if $H_\e^+\cap M=\emptyset$ for any $\e >0$, letting 
$\e \to 0$ then $M$ and $H^+$
would have a first interior point of contact, contradicting the maximum
principle. Therefore, there exists $\e >0$ such that $H_\e^+\cap
M\not=\emptyset$ and 
$(N+\e \partial/\partial t)\cap M =\emptyset$. Furthermore, the finite and
asymptotic
boundary of $H_\e^+$ is far away from $M\cup \pain M$. Consider the
hyperbolic translations along the geodesic $\beta$ with asymptotic boundary
$\{-1,1\}$, going from 1 to $-1$. Thus we would obtain a last interior
contact
point of $M$ and some translated copy of $H_\e^+$, which contradicts the maximum
principle.

 We show in the same way that $N\cap M=\emptyset$ and 
$H^-\cap M\not=\emptyset$ is not possible. This proves the Claim

\bigskip

Up to a hyperbolic translation along the geodesic $\beta$ (with asymptotic
boundary $\{-1,1\}$), we can assume that $N\cap M=\emptyset$ 
(since  $-1\in \pain \hd$ is not in the asymptotic boundary of $M$), and
therefore the Claim shows that 
$H\cap M=\emptyset$. Consider now the translated copies of $H$, along
$\beta$, going from $-1$ to 1. As $M$ is properly immersed, some translated copy
of
$H$ will have a first contact point with 
$M$ at a point $p\in M$, which
contradicts the maximum principle. This concludes the proof of the Corollary 
\end{proof}

The following result is a direct consequence of the proof of the Asymptotic
Theorem (\thmref{T.asymp.boundary}).

\begin{corollary}
Let $S_\infty \subset \pain \hd \times \R$ be a closed set strictly contained 
in a slab with width equal to $\pi$. Assume that the vertical projection of 
$S_\infty$ on $\pain \hd \times \{0\}$ omits an open arc. Then, there is
no complete properly immersed minimal surface $M$ in $\hd \times \R$ with
asymptotic boundary $S_\infty $.

\end{corollary}

\begin{proof}
By assumption there exists a complete geodesic  $c\subset \hd \times \{0\}$
such that $S_\infty$ is contained in the asymptotic boundary  of a component
of $\hd \times \R\setminus c\times \R$. We call $U$ the other component. Let 
$C$ be a catenoid, observe that any compact part of $C$ may be mapped into
$U$ by an ambiant isometry. In this situation, we can proceed as in the proof
of 
\thmref{T.asymp.boundary}.
\end{proof}

\begin{remark}{\em 
We will see in Proposition \ref{P.hyperbolic translation} that for any $t_0>\pi$
there exists a Jordan curve
$\Gamma_\infty \subset \pain \hd \times [0,t_0]$, homologous to zero, which is
the asymptotic
boundary of a properly embedded complete minimal surface. Therefore the results
in Theorem \ref{T.asymp.boundary} and \coref{C.asym.bound}
are sharp. The formulae of the generated
curves in Proposition \ref{P.hyperbolic translation} are the same as
formulae given by the first author in  \cite{Sa}. The geometric description
of the surfaces given in Proposition \ref{P.hyperbolic translation} is new.
We remark that L. Hauswirth \cite{Haus} 
has given a classification of minimal surfaces
invariant by hyperbolic translations using another approach.
}
\end{remark}

\begin{proposition}\label{P.hyperbolic translation}
Let $q_1,q_2\in \pain \hd$ $($identified with $ \pain \hd \times \{0\})$, be
two
distinct asymptotic points. Let $\ga \subset \hd$ $($identified with 
$\hd \times \{0\})$, be the complete geodesic with asymptotic boundary
$\{q_1,q_2\}$. Let us call $c_1$ $($resp. $c_2)$ the closed arc in 
$\pain \hd$ joining $q_1$ to $q_2$ $($resp. $q_2$ to $q_1)$ with respect to the
counterclockwise orientation.

  There exist a one-parameter family $M_d,\ d>0$, of complete properly embedded
minimal surfaces, invariant by the hyperbolic translations along $\ga$. The
geometric behaviour of $M_d$ is as follows.
\begin{enumerate}
\item If $d>1$, then $M_d$ contains the equidistant line $\ga_d$ of $\ga$ in 
$\hd \times \{0\}$ staying at the distance $\cosh^{-1} (d)$ from $\ga$ in the
connected component of $\hd \backslash \ga$ whose asymptotic boundary is $c_1$.
Furthermore $M_d$ is symmetric with respect to the slice $\hd \times \{0\}$ and
we have $($see Figure $\ref{F.homologous.zero})$
%$($see Figure $\ref{F.homol.zero})$
\begin{align*}
 \pain M_d &=  (c_1\times \{-H(d)\})\cup  (c_1\times \{H(d)\})  \\
  & \quad \cup\  (\{q_2\}\times [-H(d),H(d)]) \cup (\{q_1\}\times [-H(d),H(d)]),
\end{align*}
where
\begin{equation}\label{Eq.function.H}
H(d):=\int_{\cosh^{-1}( d)}^{+\infty}\dfrac{d}{\sqrt{\ch^2u-d^2}}du,\ \  d>1.
\end{equation}
%\begin{figure}[ht]
%\includegraphics[scale=0.4]{hiperbolictranslation_1.eps}
%\caption{}\label{F.homol.zero}
%\end{figure}
Therefore, $\pain M_d$ is a Jordan curve homologous to zero in $\pain\hd \times
\R$. Furthermore, the part $M_d\cap (\hd \times[0,H(d)]$ is a graph over the
component of $\hd \backslash \ga_d$ whose asymptotic boundary is $c_1$.
Finally, $H(d)$ is a nonincreasing function satisfying 
$$
\lim_{d\to 1}H(d)=+\infty,\ \ \ \lim_{d \to +\infty}H(d)=\frac{\pi}{2}.
$$

\begin{figure}[h]
\centerline{
\subfigure[]{
\includegraphics[scale=0.25]{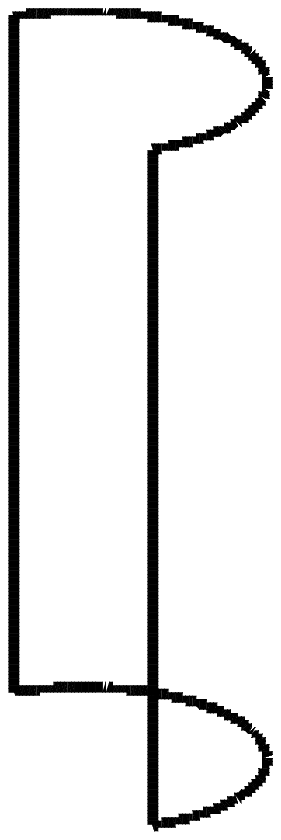}
\label{F.homologous.zero}
}
\subfigure[]{
\includegraphics[scale=0.25]{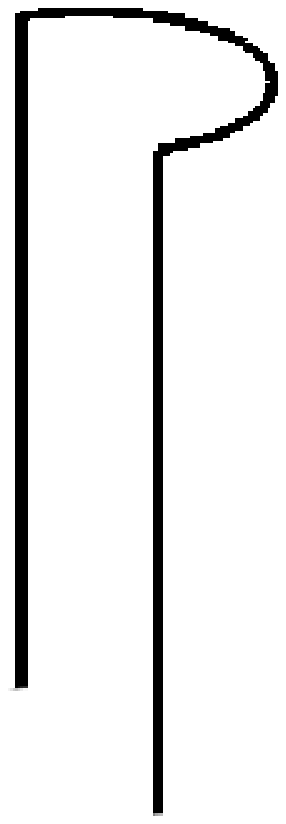}
\label{F.open.curve}
}
\subfigure[]{
\includegraphics[scale=0.25]{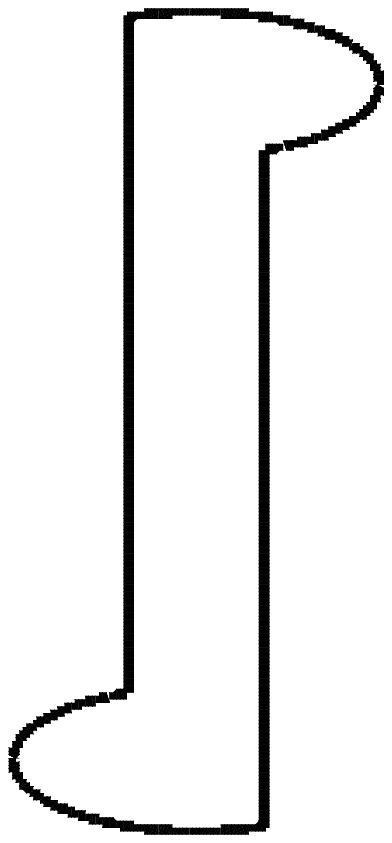}
\label{F.nonhomologous.zero}
}
}
\caption{}\label{F.hyperbolic.asymp}
\end{figure}

\item If $d=1$, then $M_1$ is the surface given by Formula $\ref{sch}$ , and its
asymptotic boundary is given by $($see Figure
$\ref{F.open.curve})$
%$($see Figure $\ref{F.open.curve})$
$$
\pain M_1= (\{q_1\} \times (-\infty,0])\cup c_1 
\cup (\{q_2\} \times (-\infty,0]).
$$
%\begin{figure}[ht]
%\includegraphics[totalheight=0.25\textheight]
%{hiperbolictranslation_3.eps}
%\caption{}\label{F.open.curve}
%\end{figure}
\item If $0<d<1$, then $M_d$ is an entire vertical graph over $\hd$ and
contains the geodesic $\ga \times \{0\}$. The asymptotic boundary of $M_d$ is
given by $($see Figure $\ref{F.nonhomologous.zero})$
\begin{align*}
 \pain M_d &= (c_2 \times \{-G(d)\}) \cup (c_1 \times \{G(d)\})  \\
  & \quad \cup \ (\{q_1\} \times [-G(d),G(d)]) \cup (\{q_2\} \times
[-G(d),G(d)]),
\end{align*}
where 
$$
G(d):=\int_0^{+\infty} \dfrac{d}{\sqrt{\ch^2u-d^2}}du,\ \ 0<d<1.
$$
%\begin{figure}[ht]
%\includegraphics[totalheight=0.23\textheight]
%{hiperbolictranslation_2.eps}
%\caption{}\label{F.nonhomol.zero}
%\end{figure}
Therefore, $\pain M_d$ is a Jordan curve non homologous to zero in 
$\pain\hd \times \R$. Furthermore $G(d)$ is a nondecreasing function and we have
$$
\lim_{d \to 0}G(d)=0,\ \ \ \lim_{d \to 1}G(d)= +\infty.
$$
\end{enumerate}
\end{proposition}

\begin{proof}

 We work with the disk model for $\hd $, so that
$$
\hd =\{(x,y)\in \m 2,\ x^2+y^2<1 \}.
$$
Therefore the product metric on $\hd \times \R$ reads as follows
$$
d \tilde s ^2=\left( \dfrac{2}{1-(x^2+y^2)}\right)^2 (d x^2
+ d y^2) + d t^2,
$$
where $(x,y)\in \hd$ and $t\in \R$.

Up to an isometry, we can assume that $q_1=-i$ and $q_2=i$. Therefore we have
$c_1=\{e^{i\theta};\ \ -\pi/2 \leqs \theta \leqs \pi/2\}$ and 
$c_2=\{e^{i\theta};\ \ \pi/2 \leqs \theta \leqs 3\pi/2\}$.

 We consider the following
particular geodesic of $\hd$
\begin{equation*}
\Gamma \!= \!\{(x,0), \ x\in (-1,1)\,\}\!\subset \! \hd .
\end{equation*}

We can assume that the surfaces invariant under
hyperbolic
translation along $\ga$ (called {\em hyperbolic surfaces}), are
generated by curves in the vertical geodesic plane $P= \Gamma
\times \R \subset \hd\times\R$.

On the geodesic $\Gamma$ we denote by $\rho \in \R$ the signed
distance to the origin $(0,0)$, thus $x=\tanh (\rho/2)$. Therefore
the metric on $P$ is
$$
d s^2=d\rho^2 +d t^2.
$$
Let us consider a curve in $P$ which is a vertical graph:
$c(\rho)=(\rho, \lambda (\rho))$ where $\lambda$ is a smooth real
function defined on a part of $\rho \geq 0$. Let us call $M$ the hyperbolic
surface generated by $c$. On $M$ we consider the orientation given by the
upward unit
normal field. With respect to this orientation the principal
curvatures of $M$ are given by
\begin{equation*}
k_1(\rho)=\frac{\lambda^{\prime \prime}}{(1+\lambda^{\prime 2})^{3/2}} (\rho),\
\ 
\mathrm{and}\ \ 
k_2(\rho)=\frac{\lambda^{\prime}}{\sqrt{1+\lambda^{\prime 2}}} (\rho) \tanh
(\rho).
\end{equation*}
So that, $M$ is a minimal surface if and only if 
\begin{equation*}
 \lambda^{\prime 2}= 
\frac{d^2}{\cosh^2 \rho-d^2},
\end{equation*}
for some $d\geqs 0$.

Up to the isometry $(z,t)\to (z,-t)$, we can assume that $\lambda$ is a
nondecreasing function, that is, $\lambda^\prime \geq 0$. Therefore, the
condition for
$M$ being minimal is 
\begin{equation}\label{Equation1}
 \lambda^\prime (\rho)=\frac{d}{\sqrt{\cosh^2 \rho -d^2}}
\end{equation}

In the case where $d>1$ we can choose, up to a vertical translation,
$$
\lambda(\rho)=\int_{\cosh^{-1} (d)}^\rho \frac{d}{\sqrt{\cosh^2 u -d^2}}du,
$$
for $\rho\geqs \cosh^{-1} (d)$. Setting $v=\cosh u/d -1$ we obtain:
$$
H(d)=\int_0^{+\infty}\frac{dv}{\sqrt{(v+1)^2-1}\sqrt{(v+1)^2-1/d^2}}.
$$
This shows that $H(d)$ is a nonincreasing function. Furthermore

\begin{align*}
 \lim_{d \to +\infty}H(d) &= 
\int_0^{+\infty}\frac{dv}{\sqrt{(v+1)^2-1}\sqrt{(v+1)^2}},\\
    &= 
\int_0^1 \frac{dx}{\sqrt{1-x^2}},\ \ \mathrm{setting}\ x=\dfrac{1}{v+1},\\
  &= \dfrac{\pi}{2}.
\end{align*}

{\em Claim.} We have $\lim_{d\to 1} \lambda (\rho)=+\infty$ for 
any $\rho >0$. 

\medskip
This clearly implies that 
$\lim_{d\to 1} H(d)=+\infty$, since 
$H(d)=\lim_{\rho\to +\infty} \lambda (\rho)$ and $\lambda (\rho)$ is a
nondecreasing function.

Setting $v=\cosh (u)/d -1$ we get
\begin{align*}
\lambda(\rho)  &= 
\int _0^{\cosh (\rho) /d -1}
\frac{1}
{\sqrt{v+2}\sqrt{(v+1+\frac{1}{d})(v+1-\frac{1}{d})}} \frac{dv}{\sqrt{v}} \\ 
  &\geqs \frac{1}{\sqrt{\frac{\cosh (\rho)}{d}+1} \sqrt{\frac{\cosh (\rho)
+1}{d}}}
\int _0^{\cosh (\rho) /d -1}\frac{dv}{\sqrt{v^2 +(1-\frac{1}{d})v}}.
\end{align*}
Denoting by $I(d)$ the last integral, we have
\begin{equation*}
I(d)= \frac{2}{1-\frac{1}{d}}
\int _0^{\cosh (\rho) /d -1} \frac{dv}
{\sqrt{( \frac{2v}{1-\frac{1}{d}}     +1    )^2 -1   }} .
\end{equation*}
Setting $s=\frac{2v}{1-\frac{1}{d}}+1$ we obtain
\begin{equation*}
I(d)=\int_1^{\frac{2\cosh (\rho) -d -1}{d-1} } \frac{ds}{\sqrt{s^2-1}} =
\cosh^{-1}( \frac{2\cosh (\rho) -d -1}{d-1}),
\end{equation*}
from what we deduce that $\lim_{d\to 1}I(d)=+\infty$ for any $\rho >0$, 
which concludes the proof of the Claim.

\medskip

In the case where $d=1$ we can choose, up to a vertical translation,
$$
\lambda(\rho)=\log (\frac{e^\rho -1}{e^\rho +1}),
$$
for $\rho >0$. The hyperbolic surface generated by $\lambda$ is a vertical
graph over the connected component of $\hd\backslash \ga$ whose asymptotic
boundary is $c_1$. This graph takes value $-\infty$ on $\ga$ and  value zero on
$c_1$ (because $\lim_{\rho \to 0}\lambda(\rho)=-\infty$ and 
$\lim_{\rho \to +\infty}\lambda(\rho)=0$). 
Since this is the unique hyperbolic surface, up to isometry, with unbounded
height, we deduce that $M_1$ is congruent to the hyperbolic surface
given by Formula  (\ref{sch}).

\medskip

Finally, in the case where $0<d<1$, the function $\lambda$ is defined for any 
$\rho \geqs 0$ and, up to a vertical translation, we can set
$$
\lambda(\rho)=\int_0^\rho \frac{d}{\sqrt{\cosh^2 u -d^2}}du.
$$
We can extend $\lambda$ on $\R$ setting 
$\lambda (\rho):=-\lambda(-\rho)$ for any $\rho\leqs 0$. Therefore $\lambda$ is
defined on $\R$ and is an odd function, and the hyperbolic surface $M_d$
generated by $\lambda$ is an entire vertical graph on $\hd$ symmetric with
respect to $\ga$. We can prove in the same way as in the case where $d>1$ 
that for any $\rho >0$ we have 
$\lim_{d \to 1}\lambda (\rho)=+\infty$. This implies that 
$\lim_{d \to 1}G(d)= +\infty$.

The other assertions in the Statement are straightforward
verifications.

\end{proof}

\section{ minimal vertical  graphs}\label{ghr}

There are many notions of graphs in $\hip^2\times \R$, but the
notion of {\em minimal vertical graphs} has appeared in   many
important theorems. See, for instance ~\cite{C-R}, ~\cite{F-M1},
~\cite{hst}, ~\cite{HRS},~\cite{Sa}, ~\cite{Sp}.

Consider  a $C^2$ function $t =u(x, y).$ The {\em vertical minimal
 equation} in $\hiper^2\times \Rea,$ is given by the following
equation:

\begin{equation}
\label{graphip} \diver_{\h}\left(\frac{\nabla_{\h}
u}{W_u}\right)=0
\end{equation}

where $\diver_{\h}$ and $\nabla_{\h} $ are the hyperbolic
divergence and gradient respectively and
$W_u=\sqrt{1+|\nabla_{\h}u|_{\h}^2},$ being  $|\cdot|_{\h}$ the
norm in  $\h^2.$

Focusing the halfspace model for $\h^2,$ with Euclidean
coordinates $x,y,$  $y>0,$ the {\em minimal vertical  equation}
\eqref{graphip} takes the following form
\begin{align}\label{mea2}
\begin{split}
(1 + y^2u_x^2)u_{yy} + (1 + y^2u_y^2)u_{xx} -2 y^2u_x u_y u_{xy}
-yu_y(u_x^2 + u_y^2)=0
\end{split}
\end{align}

There are many explicit examples of entire and complete minimal graphs with
nice geometric
properties. For instance,  in the half-plane
model,

\begin{enumerate}

\item The equation \cite{Sa}
\begin{equation*}\label{csch} t = \ell
x,\, x\in (-\infty, \infty),\  y>0
\end{equation*}
 gives rise to an entire minimal graph (left side of Figure \ref{F.examples})
symmetric 
about the geodesic $\{x=0\}$, that is constant (the
constant varying in the interval $(-\infty, \infty)$)
 on each leaf of the foliation  given by geodesics with a fixed common
asymptotic
boundary point $p$ (in this~ model~$p=\infty$). Thus the
asymptotic boundary consists in the union of a vertical line with
a complete embedded curve  in $\pain \hip^2\times \R$ asymptotic
to that line.
%\begin{figure}[ht]
%\includegraphics[scale=0.45]{m11111.eps}
%\caption{ Ball model for $\hipe^2\times\{0\}$}
%\end{figure}
%\begin{multline*}
%\epsfxsize=7cm \epsfbox{m11111.eps}\\
%\hfill\text{\bf Figure 1: ball model for $\hipe^2\times\{0\}$}\hfill\\
%\end{multline*}
%\includegraphics<1->[height=5.5cm]{m11111.pdf}\hfill

\item The equation \cite{Sa}
\begin{equation*}\label{log}
t=\frac{\ell}{2} \ln(x^2 +y^2),\qquad y>0
\end{equation*}
yields  an entire minimal graph (right side of Figure \ref{F.examples})
symmetric about the
geodesic 
$\{ x^2 +y^2=1,\ y>0\}$, 
that is constant on each leaf of the
foliation given by (hyperbolic) translations of a fixed geodesic;
hence, the asymptotic boundary consists of two
 embedded curves in $\pain \hip^2\times \R$ with two symmetric
ends, each end asymptotic to a half-vertical line.

%\begin{figure}[htp]
%\begin{center}
%\subfigure[]{\label{Parabolic}
%\includegraphics[scale=0.45]{m11111.eps}
%\subfigure[]{\label{Hyperbolic}
%\includegraphics[scale=0.45]{min2.eps}
%\end{center}
%\caption{Ball model for $\hipe^2\times\{0\}$}\label{F.examples}
%\end{figure}

\begin{figure}[ht]
%\begin{center}
\includegraphics[scale=0.40]{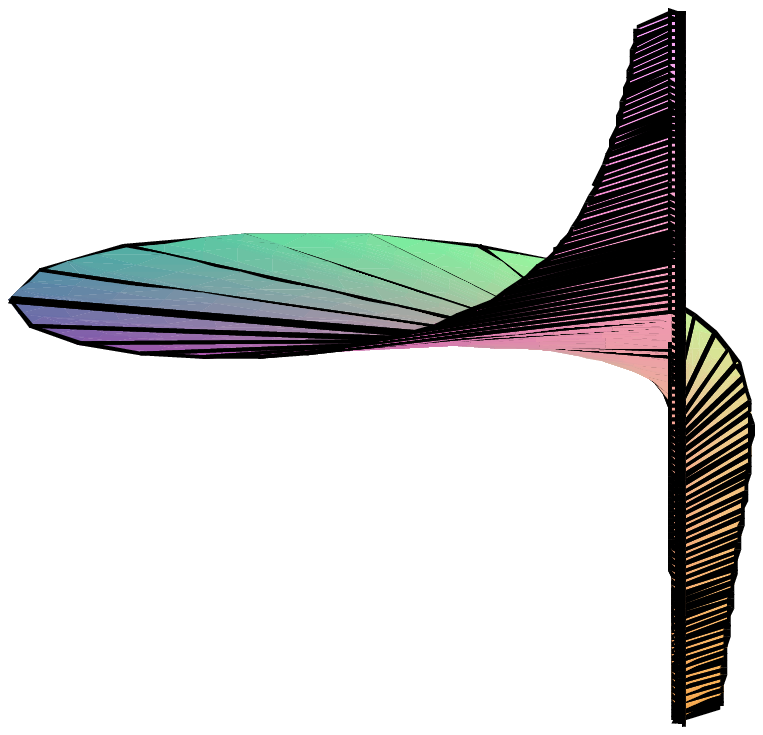}
\includegraphics[scale=0.40]{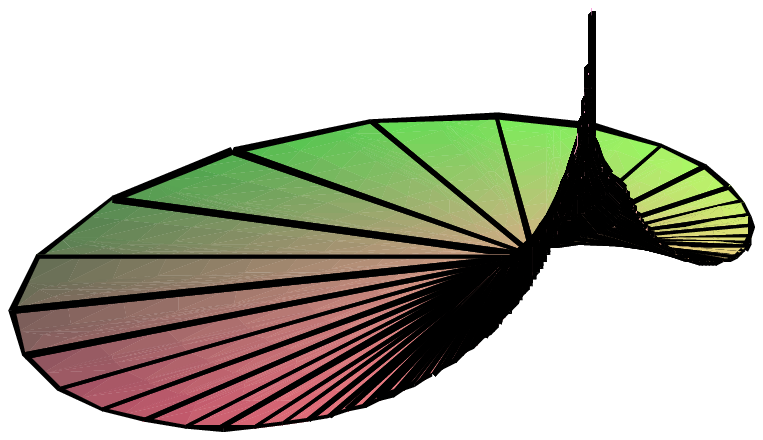}
%\end{center}
\caption{Ball model for $\hipe^2\times\{0\}$}\label{F.examples}
\end{figure}

%\begin{multline*}
%\epsfxsize=7cm \epsfbox{min2.eps}\\
%\hfill\text{\bf Figure 2: ball model for $\hipe^2\times\{0\}$}\hfill\\
%\end{multline*}

%\includegraphics<1->[height=5.5cm]{min2.pdf}\hfill

 Of course, the previous examples give two different explicit non
trivial minimal graphs over a half-plane of $\hd$ taking zero boundary
value data on a geodesic but having different asymptotic
boundaries.

\item We observe that there exists a function which takes infinite
boundary value data on the positive $y$ axis and zero asymptotic
value boundary data at the positive $x$ axis (halfspace model for
$\h^2$), invariant by {\em hyperbolic translations} \cite{Sa}.
\begin{equation}\label{sch}
t = \ln \left( \frac{\sqrt{x^2 +y^2} + y}{x}\right),\qquad y>0,
x>0
\end{equation}

\end{enumerate}
Notice that \eqref{sch} yields a complete vertical minimal graph
over a
 domain bounded  by a geodesic in $\hip^2\times \{0\},$ taking infinite
boundary value data on the geodesic  and zero asymptotic boundary
value
 data on an arc $L$ of $\partial_\infty\hip^2\times
\{0\}$. The asymptotic boundary  of the graph is then the union of
$L$ with the two upper half vertical lines arising from the end points
of $L.$ We will use this special minimal vertical graph as a {\em
barrier} at an asymptotic boundary point.

We remark that the miniaml surface given by
Formula \eqref{sch}   was used by P. Collin and H.
 Rosenberg  \cite{C-R} in the important construction of entire minimal graphs in
$\hip^2\times \R$ that are conformally the complex plane $\Cx$,
  disproving a conjecture by R. Schoen.

In the Poincaré disk model for $\hip^2$ with Euclidean coordinates
$x, y,$ $x^2 +y^2 <1,$ the {\em minimal vertical equation}
\eqref{graphip} becomes
\begin{multline}\label{meaq2}
 \Dc (u):=( 1 + \frac{(1 -x^2 -y^2)^2}{4}\, u_x^2) u_{yy}+ 
(1 + \frac{(1 -x^2 -y^2)^2}{4}\, u_y^2) u_{xx}\\
  - 2\, \frac{(1
-x^2 -y^2)^2}{4}\, u_x u_y\, u_{xy} + 2\frac{(1 -x^2
-y^2)}{4}\left (xu_x +yu_y\right) (u_x^2 +u_y^2)=0\hfill
\end{multline}
We observe that equation \eqref{meaq2} is a second order
quasilinear strictly elliptic equation for all {\em real values}
of  the independent variables $x, y$. Moreover, the eigenvalue of
the associate matrix are $1$ and\\ $W_u=1+\frac{(1 -x^2
-y^2)^2}{4} (u_x^2 +u_y^2).$ The same observation holds for the
equation
\eqref{mea2}, replacing $\frac{(1 -x^2 -y^2)^2}{4}$ by $y^2.$
Hence we conclude that both equations are regular and strictly
(uniformly) elliptic up to the asymptotic boundary of
$\hip^2\times \{0\}.$ For this reason we can state the classical
maximum principle and uniqueness for prescribed continuous finite
and asymptotic boundary data.

\begin{theorem}[Classical maximum principle]\label{mp}
Let $g_1, g_2:\ \partial \Omega \cup \pain \Omega \rightarrow \R$ be continuous
functions satisfying $g_1 \leqs g_2$. Let 
$u_i : \overline \Omega \rightarrow \R$ be a continuous extension of $g_i$
satisfying the minimal equation $(\ref{meaq2})$ on $\Omega$, $i=1,2$. Then
$u_1\leqs u_2$.

\end{theorem}

\begin{proof}
The proof is classical elliptic theory, since the minimal equation
(\ref{meaq2}) is strictly elliptic up to the asymptotic boundary.
 A geometric approach can
be done in this way. Assume that $u_1(p) >u_2(p)$ at some point 
$p\in \Omega$. Then, lifting the graph of $u_2$ vertically we obtain a last
interior contact point between the graph of $u_1$ and the graph of $u_2$, which
gives a contradiction by the interior maximum principle. 

\end{proof}

 We will  solve some Dirichlet problems over certain
unbounded domains, given certain prescribed finite boundary data
and given certain prescribed asymptotic boundary data.

 Among such domains we will consider  exterior domains $\Om$.   Of course, the classical examples of such minimal
  graphs over
  an exterior domain are given by the one parameter family of
  half-catenoids, see \lemref{l1}. We will use this family as {\em barriers}. We show some generating curves
  in Figure 3, where $R=\tanh \rho/2,$ and $\rho$ is the
  hyperbolic distance from the axe $t.$

%\begin{multline*}
%\epsfxsize=7cm \epsfbox{catenoids2_R.eps}\\
%\hfill\text{\bf Figure 3: ball model for $\hipe^2\times\{0\}$}\hfill\\
%\end{multline*}

\begin{figure}[ht]
\includegraphics[scale=0.5]{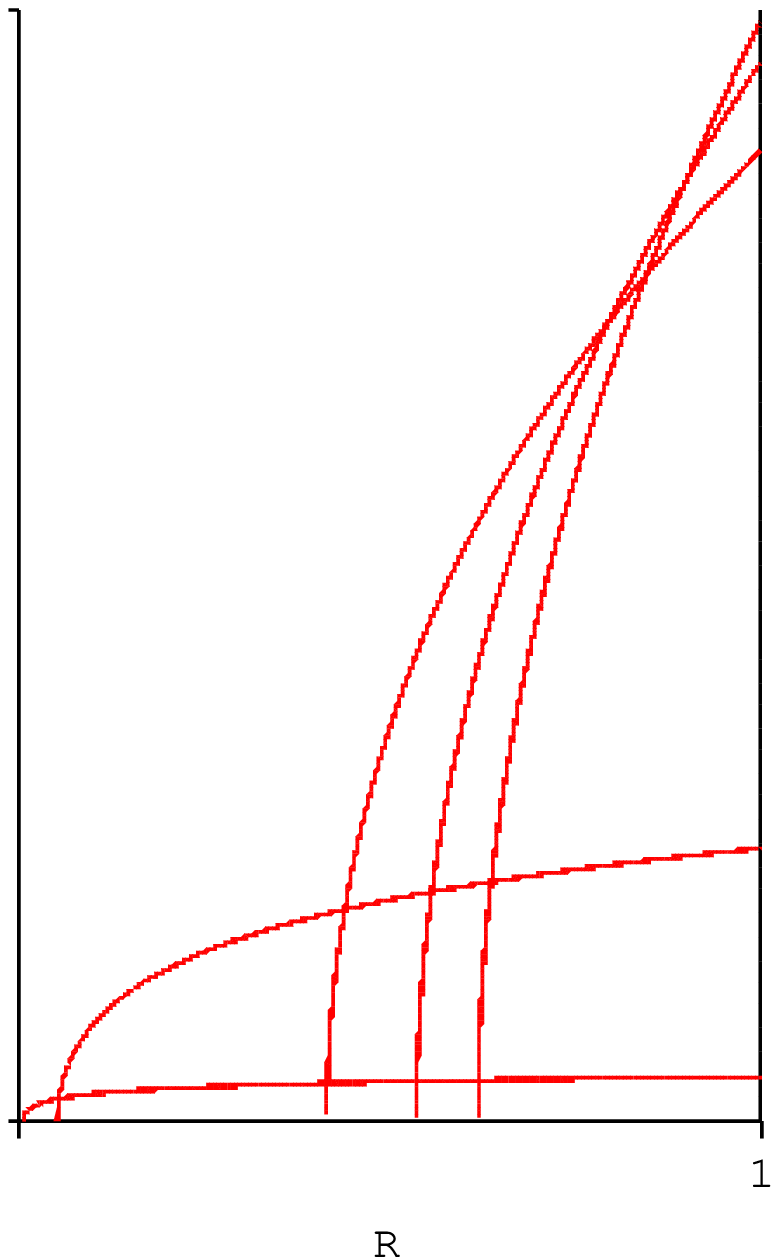}
\caption{ Ball model for $\hipe^2\times\{0\}$}
\end{figure}

 We remark that we use also as {\em barriers} the one-parameter
 family of minimal surfaces invariant by hyperbolic translations
given by Proposition \ref{P.hyperbolic translation}.

\section{The Perron process for the minimal vertical
equation}\label{s2}

In the product $\hip^2\times \R,$ we consider the disk model for
the hyperbolic plane $\hip^2$. Let $\Om\subset \hip^2\times
\{0\},$ be a domain. In $\overl{\hip^2}\times \{0\},$ we have that
$\partial \overl{\Om}=\partial \Om\cup\partial_\infty\Om,$ where
$\partial \Om\subset \hip^2\times \{0\}$ and
$\partial_\infty\Om\subset\partial_\infty \hip^2\times \{0\}.$

\begin{definition}[Problem ($P$)]\label{probP}{\em  Let 
$g:\partial \Om\cup\partial_\infty\Om \goto \R$ be a
continuous function except maybe at a finite set $S$   of  points
(discontinuities). We consider the Dirichlet problem, say Problem
$(P)$, for the minimal vertical equation \eqref{meaq2} taking at
any point of $\partial\Om\cup\partial_\infty \Om\setminus S$,
prescribed boundary (finite and asymptotic)  value data $g.$
}
\end{definition}

Let $u: \overline {\Omega}_S:= \overline \Omega\setminus S
\rightarrow \R$ be a continuous function.

Let $U\subset \Omega$
be a closed round disk in $\hip^2\times\{0\}$. If 
$u_{ |\partial U}$ is a $C^1$ function then solving the Plateau problem
~\cite{M}
and using a standard adaptation of Rado's Theorem ~\cite{R} (since
$u_{|\partial U}$ is a vertical graph over a circle), it follows
that $u_{|\partial U}$ has an unique minimal 
extension $\tilde u$ on $U$, continuous up to 
$\partial U$. If $u_{|\partial U}$
is $C^0$,
one uses an
approximation argument or uses a local barrier at a boundary point
of $U.$  We then define the continuous function $M_{U}(u)$ on
$\overline{\Omega}_S$ by:
\begin{align}\label{mu}
M_{U}(u)(x)=
\begin{cases}
      u(x) &\,\,\, \text{ if $x \in \overline{\Omega}_S\setminus U$}\\
   \tilde u (x)& \,\,\, \text{ if $x \in U$}
  \end{cases}
\end{align}
%\vskip2mm

We say that $u$ is a {\em subsolution} (resp. {\em
supersolution})
of $(P)$ if:
\begin{itemize}
\item [i)]\label{Item.round.disk} For any closed round disk
$U \subset \Omega$ we have\\ $u \leqs M_{U}(u)$ (resp. $u \geqs
M_{U}(u)$).  \vskip1mm
\item [ii)]  $u\mid_{\partial\Om\cup \partial_\infty\Om}\; \leqs g$
(resp. $u\mid_{\partial\Om\cup \partial_\infty\Om}\;\geqs g)$.

\end{itemize}

\begin{remark}\label{r1}{\em
 We now give some classical facts about subsolutions and
supersolutions, see \cite{C-H}, \cite{asi}.
\begin{enumerate}
 \item  It is easily seen that if $u$ is $C^2$ on $\Omega$, the condition 
{ i})  
above  is equivalent to $\dc u
\geqs 0$ for subsolution or $\dc u \leq 0$ for supersolution.
 \item As usual if $u$ and $v$ are two subsolutions (resp.
supersolutions) of $(P)$ then $\sup (u,v)$ (resp. $\inf (u,v)$)
again is a subsolution (resp. supersolution).
 \item Also if $u$ is a subsolution (resp. supersolution) and
$U \subset \Omega$ is a closed round disk then $M_{U}(u)$ is again
a subsolution (resp. supersolution).

\item \label{Item.bounded} Let $\phi $  (resp. $u$) be a supersolution  (resp. a
subsolution) of
Problem $(P)$ such that $u\leqs \phi$, then we have 
$M_U(u)\leqs M_U(\phi)\leqs \phi$ for any disk 
$U$ with $\overline U \subset \Omega$.
 
Note that if 
$\phi$ and $u$ are continuous on $\overline \Omega$ then
necessarily $u\leqs \phi$ on $\Omega$.

\end{enumerate}
}
\end{remark}

Note also that 
due to the nature of Equation \eqref{meaq2}, $\Omega$ is a bounded domain in 
$\overline {\hd}\times \{0\}$. 

\begin{definition}[Barriers]\label{d1}{\em

We consider the Dirichlet Problem $(P)$, see \defref{probP}. Let
$p \in \partial \Om\cup\partial_\infty\Om, $ be a boundary point where
$g$ is continuous.

\begin{enumerate}
\item \label{D.Item.bar.1}
 Suppose that for any $M>0$ and 
for any $k\in\Na$ 
there is an open neighborhood $ \Nc_k$ of $p$ in
$\R^2$ and a function
$\omega_{k}^+$ $($resp. $\omega_{k}^-)$  in 
$C^{2}(\Nc_k \cap \Omega)\cap C^{0}(\overline{\Nc_k \cap \Omega)}$  
such that
\begin{itemize}
\item[ i)]$\omega_{k}^+(x)\mid_{(\PO) \cap \Nc_k}\geq g(x)$ and
$\omega_{k}^+(x)\mid_{\partial \Nc_k  \cap \Omega} \geq M$\\
$($resp. $\omega_{k}^-(x)\mid_{(\PO) \cap \Nc_k}\leq g(x)$ and
$\omega_{k}^-(x)\mid_{\partial \Nc_k  \cap \Omega} \leq -M)$
\item[ii)] $\Dc (\omega_{k}^+) \leq 0$ $($resp.
$\Dc (\omega_{k}^-) \geq 0)$
 in $\Nc_k \cap \Omega$,
\item[ iii)] $\lim_{k \to +\infty}\omega_{k}^+(p)=g(p)$
 $($resp. $\lim_{k \to +\infty}\omega_{k}^-(p)=g(p))$.
\end{itemize}

\item \label{D.Item.bar.2} Suppose that there exists a supersolution $\phi$
$($resp.
a subsolution $\eta)$ in $C^2(\Om)\cap C^0(\overline \Om)$ such that
$\phi(p)=g(p)$ (resp. $\vphi(p)=g(p)).$

\end{enumerate}
  In both cases (\ref{D.Item.bar.1}) or  (\ref{D.Item.bar.2})
we say that $p$ admits a
{\em superior barrier} ($\omega_{k}^+, k\in \Na$  or $\phi$)
(resp. {\em inferior} {\em barrier} $\omega_{k}^- , k\in \Na $ or
$\vphi$) for the Problem ($P$). If $p$ admits a superior and an
inferior barrier we say more shortly that $p$ admits a {\em
barrier}.}

\end{definition}

\begin{example}{\em [Barrier
at any convex point  for any   bounded continuous boundary data
$g$]}. \label{bf}{\em The construction of B. Nelli and H.
Rosenberg, the Scherk type minimal graph in $\hip^2\times\R$ over
a geodesic triangle, taking zero values data on two sides and
infinite at the other side ~\cite{nel} is given in
\exaref{scherk}. The geodesic triangle and the boundary data are
drawn in Figure 4.

\begin{figure}[ht]
\includegraphics[scale=0.5]{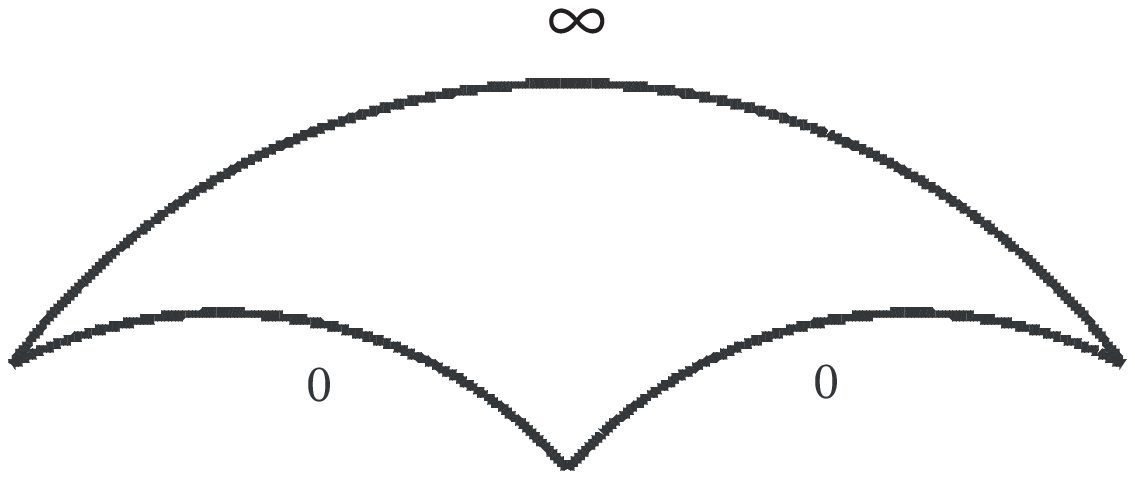}
\caption{}
\end{figure}

%\begin{multline*}
%\hfill\epsfxsize=5cm \epsfbox{scherk_NR.eps}\hfill\\
%\hfill \text{\bf Figure 4}\hfill
%\end{multline*}

We consider now these Scherk type surfaces when the geodesic
triangle $\Delta$ is isosceles and the zero value data are taken
on the two sides with equal length and $-\infty$ on the other
side. We show that these surfaces can be used as a upper barrier
(in the sense of \defref{d1}-(\ref{D.Item.bar.1}))
at any convex point $p_0\in\partial \Om,$ for any  boundary
bounded data $g$ continuous at $p_0.$
 For the lower barrier the construction is analogous. 

Let
$\Delta$ be a geodesic isosceles triangle in $\hip^2\times\{0\}$
with sides $A$, $C_1$ and $C_2,$ with $|C_1|= |C_2|.$ Let $\om$ be
the solution of the minimal equation
 taking zero value data on $C_1$ and $C_2,$ and
$-\infty$ on $A.$  Let $S$ be the graph of $\om$. Let $a$ be the
common vertex of $C_1$ and $C_2.$ Let $\ga $ be the axe of
symmetry of $\De,$ hence $a\in \ga.$ Let $\{b\}=\ga\cap A.$ 
Let $\be$ be a
geodesic intersecting $\Delta$  orthogonal to $\ga$ at a point
$d\in (a, b)$. Set  $\be\cap C_1=\{c\}.$

We claim the following:
\begin{enumerate}

\item $\om$ along $\ga$ is nonincreasing in $[a, b].$

\item $\om$ along $\ga$ is nonincreasing in $[c, d].$

\end{enumerate}

Assume momentarily the Claim.  Let $p_0\in\partial \Om$ be a
convex point and let $g$ be a boundary data continuous at $p_0$.
Let $M>0$ be any positive real number. It suffices  to show that
for any $k\in\Na$ there is an open neighborhood $ \Nc_k$ of $p_0$
in $\R^2$   and a function $\omega_{k}^+$   in $C^{2}(\Nc_k \cap
\Omega)\cap C^{0}(\overline{\Nc_k \cap \Omega)}$ such that
\begin{itemize}
\item[ i)]$\omega_{k}^+(x)\mid_{\partial\Om \cap \Nc_k}\geq g(x)$ and
$\omega_{k}^+(x)\mid_{\partial \Nc_k  \cap \Omega} \geq M$\\

\item[ii)] $\Dc (\omega_{k}^+) = 0$
 in $\Nc_k \cap \Omega$,
\item[ iii)] $\omega_{k}^+(p_0)=g(p_0) +1/k$.
\end{itemize}

 By continuity there exists $\ep>0$ such that for any $p\in\partial
\Om$ such that $\dist(p, p_0) <\ep$ we have $g(p) <g(p_0) + 1/k.$
By assumption there exists an open geodesic arc $\ga^\perp,$
through $p_0$ such $ \ga^\perp\cap\Om=\emptyset.$ We may assume
that the disk $D_\ep(p_0)$ intersects $\ga^\perp$ at two points.

We choose $\De$ such that $p_0\in \ga$, $|A| <\ep,$ 
$\ov{\Om}\cap A=\emptyset,$ and $\ga$ orthogonal to 
 $\ga^\perp$ at $p_0.$ Let
$M_1>\max\{M, g(p_0) +1/k\}.$ We consider the Scherk surface 
(graph of $\om$) 
taking $M_1$ boundary value data on $C_1$, $C_2$ and  $-\infty$ on
$A.$ By continuity, there exists a point $p_1$ at $\ga$ where $\om
(p_1)=g(p_0) +1/k.$ Up to a horizontal translation along $\ga$
sending $p_1$ to $p_0$, we may assume that $\om(p_0)=g(p_0) +1/k.$
Therefore we set $\Nc_k=\De\cap\Om$ and $\om_k^+=\om\mid_{\Nc_k}$
is the restriction of $\om$ to $\Nc_k.$  The Claim shows that
$\omega_{k}^+(x)\mid_{\partial\Om \cap \Nc_k}\geq g(x)$, as
desired.

We will now proceed the proof of the Claim. Let $p_1, p_2\in [a,
b)$ such that $p_1 <p_2.$ Let $p_3\in (p_1, p_2)$ be the middle
point in the segment $[p_1, p_2]$ and let $\ga_3$ be the geodesic
orthogonal to $\ga$ at $p_3.$ Let $S_3$ be the connected component
of $S\setminus (\ga_3\times \R)$ containing $(a, 0).$  Now
maximum principle shows that  symmetric of $S_3$  with respect to
$\ga_3\times \R$ is above $S$, since this  is true on the
boundary.  Hence $\om(p_1)
> \om (p_2),$ as desired. The  proof of the second  part of the
Claim is analogous, considering the reflections about the vertical
geodesic
planes orthogonal to $[c, d]$. The same argument also shows  that 
$S$ is symmetric
about the vertical geodesic plane $\ga\times \R$. 
 This accomplishes the construction
of the desired barrier.

}
\end{example}

\begin{example}[Barrier at an asymptotic point]
\label{ba}

{\em  The surface given by Formula \eqref{sch}, may be seen as a
complete vertical minimal graph over a
 domain bounded  by a geodesic in $\hip^2\times \{0\},$ taking infinite
boundary value data on the geodesic  and zero asymptotic boundary
value
 data on an arc $L$ of $\partial_\infty\hip^2\times
\{0\}$. The asymptotic boundary  of the graph is then the union of
$L$ with the two half vertical lines arising from the end points
of $L$, see Figure \ref{F.open.curve}.
 We can therefore choose the geodesic as small as we wish
in the Euclidean sense, because the minimal equation extends
smoothly to $\overl{\hip^2}\times \{0\}.$ Then we can put a copy
of it above and bellow  the graph of $g$ at any point $p$ where $g$ is
continuous. Thus we obtain a barrier at any point $p$ of
$\partial_\infty\Om$ where $g$ is continuous,
in the sense of \defref{d1}-(\ref{D.Item.bar.1}).
}
\end{example}

\begin{theorem}[Perron process]\label{perron}

Let $\Omega \subset \hip^2\times\{0\}$ be a domain and let  $g:
\partial\Om\cup\partial_\infty\Om \rightarrow \R$ be a  continuous function
except, maybe, at a finite set $S$.
Suppose that the Dirichlet
Problem $(P)$ has a supersolution $\phi$. 
Set 
$\Sc_\phi = \{\vphi,\, \text{subsolution of }  \, (P),\, \, \vphi \leq \phi
\}$. 
Assume that $\Sc_\phi\not=\emptyset.$ We define for each $x \in
\overline \Omega\setminus S$
$$
 u(x)=\sup_{\vphi \in \Sc_\phi}{\vphi(x)}.
$$
We have the following:
\begin{enumerate}
\item \label{Item.mineq} The function $u$ is $C^2$ on $\Omega$ and satisfies the
minimal equation \eqref{meaq2}.
\item \label{Item.fin.bound} Let  $p \in \partial\Om$ be a finite boundary
point where $g$ is continuous. Suppose
  that $p$ admits a 
{\em barrier} in the sense of Definition $\ref{d1}-(\ref{D.Item.bar.1})$.
Then the solution $u$ is  continuous at $p$ and
satisfies $u(p)=g(p).$  In particular, if $ \partial\Om$ is convex at $p$ then
$u$ extends continuously  at $p$ and $u(p)=g(p).$

\item \label{Item.asym.bound} Let $p \in \partial_\infty\Om$ be an asymptotic
boundary
point where $g$ is continuous. Then  $p$ admits a {\em barrier},
$u$ is continuous at $p$ and satisfies $u(p)=g(p);$ that is, if
$(x_n)$ is a sequence in $\hip^2\times\{0\}$ such that $x_n\goto
p$ in the Euclidean sense then $u(x_n)\goto g(p).$ Particularly,
if $g$ is continuous on $\partial_\infty\Om$ then the asymptotic
boundary of the graph of $u$ is  the restriction of the graph of
$g$ to $\partial_\infty\Om.$

\item \label{Item.segment} Let 
$q\in \pain \hd$ be an interior point of 
$\partial_\infty\Om $ where $g$ is
discontinuous. Then the vertical segment \newline 
$\{(q, t),\; t\in [A:=\liminf\limits_{ x\goto q,\ x\not= q}g(x),
\, B:=\limsup\limits_{ x\goto q,\ x\not=q}
g(x)],\, x\in\partial\Om\cup~\partial_\infty\Om\}$\newline
 belongs
to the asymptotic boundary
of the graph of $u$. In particular, if $A=-\infty$ and $B=+\infty$, then the
whole vertical line $\{q\} \times \R$ belongs to the asymptotic boundary.

\end{enumerate}

\end{theorem}

\begin{proof}
{}
%\begin{enumerate}
 Observe that for any $\vphi\in\Sc_\phi$, $M_U(\vphi)\in
\Sc_\phi,$ for any closed disk $U\subset\Om.$ Observe also that
the basic {\em compactness theorem} holds for the minimal vertical
equation, see \cite{GT}, ~\cite{Si}, ~\cite{Sp} and ~\cite{HRS}.
The proof of Statements (\ref{Item.mineq}) and (\ref{Item.fin.bound}) 
follows from
classical arguments as in Theorem 3.4 in \cite{asi}, see also the
classical reference \cite{C-H}. The last assertion of Statement 
(\ref{Item.fin.bound}) follows from \exaref{bf}.

The Statement (\ref{Item.asym.bound}) follows from the previous construction
of a suitable barrier, in the sense of \defref{d1}-(\ref{D.Item.bar.1}),
at any point
$p$ of
$\partial_\infty\Om$ where $g$ is continuous, 
see \exaref{ba}.

The proof of Statement (\ref{Item.segment})
follows from a continuity
argument. Indeed, as $g$  is discontinuous at $q$ we have
$A\not=B.$
 let $t_0\in (A,B).$ Let $(x_n)$ and $(y_n),$ $\, n\in \Na,$
be two sequences in $\partial_\infty\Om,$ such that
$x_n, y_n \goto q$, $\lim g(x_n)= A,$ and $g(y_n)=B.$ We can
assume that $g(x_n) < t_0 < g(y_n),$ for any $n.$ Let $\Ga_n$ be a
closed arc joining in ${\overl\Om}$ the point $x_n$ to $y_n,$
close to $q$ in the Euclidean sense and such that 
$\Gamma_n \cap \pain \Om=\{x_n,y_n\}$.
 Notice that the restriction
of the graph of $u$ to  the closed arc $\Ga_n$ is continuous and
intersects the   slice $\hip^2\times \{t_0\},$ at some point 
$(z_n,t_0),$ where $z_n$ is an interior point of $\Ga_n,$ 
and $z_n \goto q$ as $n\goto \infty.$ Hence, $(q, t_0)$ belongs to the
asymptotic
boundary of the graph of $u$, for any $t_0\in [A, B]$ (as the
asymptotic boundary is a closed set).
This completes the proof of the Theorem.

\end{proof}

\begin{corollary}\label{C.bounded.data} 
 Let $\Omega \subset \hip^2\times\{0\}$ be a domain and let \newline 
$g: \partial\Om\cup\partial_\infty\Om \rightarrow \R$ be a  
bounded function everywhere continuous except maybe at a finite set
$S\subset \partial\Om\cup\partial_\infty\Om$. 
Assume that the finite boundary $\partial \Omega$
is convex or, alternatively, that each finite boundary point admits a barrier.

Then, $g$ admits an extension $u: \overline \Omega \setminus S\rightarrow \R$
satisfying the minimal vertical equation $(\ref{meaq2})$. Furthermore, 
the total boundary  of the
graph of $u$ $($that is the finite and asymptotic boundary$)$ is the union of
the
graph of $g$ on 
$(\partial\Om\cup\partial_\infty\Om) \setminus S$ with
the vertical segments \newline
$\{(q, t),\; t\in [A:=\liminf\limits_{ x\goto q,\ x\not= q}g(x),
\, B:=\limsup\limits_{ x\goto q,\ x\not=q}
g(x)],\, x\in\partial\Om\cup~\partial_\infty\Om\}$\newline
at any $q\in S$.

\end{corollary}

\begin{proof}
Since $g$ is bounded, there are some constant functions which are
supersolutions and other which are subsolutions of Problem ($P$). 
We consider a slight variation of Perron process taking the 
set $\mathcal S$ of continuous subsolutions of $(P)$. Let $u$ be the solution
given by the Perron process (Theorem \ref{perron}).
It follows from Theorem \ref{perron} that the total boundary of the graph
 of $u$ contains the union of the graph of $g$ on 
$(\partial\Om\cup\partial_\infty\Om) \setminus S$ with
the vertical segments given in the Statement at any 
$q\in S\cap \pain \Omega$. If $q\in S$ is on $\partial \Omega$ or is not an
interior point of $\pain \Omega$ then, taking into account that each finite
boundary point has a barrier by assumption, we can prove in the same way that
the vertical segment $[A,B]$ is contained in the total boundary of the graph of
$u$.

 For any $q_i\in S$ we set 
$A_i:=\liminf\limits_{ x\goto q_i,\ x\not= q_i}g(x)$ and 
$B_i:=\limsup\limits_{ x\goto q_i,\ x\not=q_i}g(x)$, 
$x\in \partial\Om\cup~\partial_\infty\Om$.

It remains to show that
for any $q_i\in S$ and any real number $t$ satisfying $t>B_i$ or 
$t<A_i$ the point $(q_i,t)$ is not in the total boundary of the graph of $u$.

Assume first that $t>B_i$. Let $\varepsilon >0$ be a real number satisfying 
$B_i +\varepsilon <t$. There exists a continuous function
$g^+ :\partial \Omega \cup \pain \Omega\rightarrow \R$ such that $g^+ >g$ on 
$(\partial \Omega \cup \pain \Omega)\setminus S$ and $g^+(q_j)=q_j+\varepsilon$ 
for any $q_j \in S$. Then the minimal extension $u^+$ of $g^+$ given by the
Perron
process is continuous up to $\overline{\Omega}$. It follows that $u^+$ is a
supersolution of Problem $(P)$ for the boundary data $g$
and, consequently we have $\varphi \leqs u^+$ on 
$\overline{\Omega}$ for any $\varphi \in \mathcal{S}$.
It follows that the point $(q_i,t)$ is not in the total boundary of the graph of
$u$.  

Assume now that  $t<A_i$ and consider a continuous function 
 $g^- :\partial \Omega \cup \pain \Omega\rightarrow \R$ such that $g^- <g$ on 
$(\partial \Omega \cup \pain \Omega)\setminus S$ and $g^-(q_j)=q_j-\varepsilon$ 
for any $q_j \in S$. Since the minimal extension of $g^-$ is a subsolution of
Problem ($P$),  we infer that 
the point $(q_i,t)$ is not in the total boundary of the graph of
$u$. This concludes the proof of the Corollary.

\end{proof}

\begin{remark}\label{remark1}\hfill

\begin{enumerate}
{\em 
 \item It follows from Corollary \ref{C.bounded.data} that 
if $\Om$ is a convex unbounded domain, then
there exists an unique minimal vertical graph over $\Om$ taking any
prescribed bounded continuous finite and asymptotic  boundary
data.

 \item  In the  special case when $\Omega=\hd$, consider 
a  bounded function $g$ on  $ \pain \hd\times\{0\}$, continuous 
except maybe at a finite set of points $S$. With the aid  of Corollary
\ref{C.bounded.data} we see that $g$ admits a minimal entire extension $u$.
 If $g$ is continuous, we remark that 
 uniqueness of the extension follows from \thmref{mp}.

This problem when $g$ is continuous on $\pain \hd$, is called 
{\em  Dirichlet problem at infinity}  and was solved by
 B. Nelli and H. Rosenberg ~\cite{nel}.}
\end{enumerate}

\end{remark}

\begin{example}\label{scherk}{\em
 Let $\Delta$ be a geodesic triangle in $\hip^2\times\{0\}$ with
sides $A$, $C_1$ and $C_2.$ We want to show that there exists a
minimal Scherk type graph over $\Delta$ taking zero boundary value
data on the interior of $C_1 \cup C_2$ and taking $+\infty$ as
boundary value data on $A.$ This is proved by B. Nelli and H. Rosenberg in
~\cite{nel}. 

For this purpose we first show that for any $n\in \Na$ 
there exists a
solution $u_n$ of the minimal equation on
the interior of $\Delta$ taking zero boundary value
data on the interior of $C_1 \cup C_2$ and taking $n$ as
boundary value data on $A$. We consider the
set $\Sc_n$ of  continuous functions $\vphi$  on
$\Delta$ satisfying:
\begin{enumerate}
\item For any closed round disk $U \subset \inter\Delta, \vphi \leqs
M_{U}(\vphi),$ where $ M_{U}(\vphi)$ is given in Formula
\eqref{mu}
\item $\vphi\leqs 0\quad \text{on  the interior of $C_1\cup C_2$}$
\item $\vphi \leqs n$ on $A$.
\end{enumerate}

For any subarc $C^\prime$ of $C_1 \cup C_2$ and any subarc $A^\prime$ of
$A$
there is continuous subsolutions and supersolutions on $\Delta$ assuming 
zero boundary value data on $C^\prime$ and $n$ 
boundary value data on $A^\prime$. Those functions give barriers at any 
interior point of the sides $A$, $C_1$ and $C_2$.  Therefore the solution
$u_n$ given by the Perron process, 
Theorem \ref{perron}-(\ref{Item.mineq}), assumes the desired boundary value
data

 Let
$A_\infty$ be the complete geodesic containing $A$. Taking into
account Formula \eqref{sch},  let $\phi$ be the minimal graph over
the half-plane with boundary $A_\infty$ that contains $\Delta$,
taking $+\infty$ as boundary value data on $A_\infty$ and zero
asymptotic boundary value data. We will write down
 a slight variation of Perron process.

 Let $\Sc_\phi$ be the family of continuous functions $\vphi$ defined
 on\\
 $\inter\Delta\cup \inter(C_1\cup C_2)$ satisfying:
\begin{enumerate}
\item For any closed round disk $U \subset \inter\Delta, \vphi \leqs
M_{U}(\vphi),$ where $ M_{U}(\vphi)$ is given in Formula
\eqref{mu}
\item $\vphi\leqs 0\quad \text{on  the interior of $C_1\cup C_2$}$
\item $\vphi
\leq \phi$
\end{enumerate}

 Notice that  the functions $u_n$ construted above
belong to $\Sc_\phi.$ Therefore we infer that the solution $u$ given
by Perron process assumes infinite boundary value data on $A.$ We
claim that $u$ takes zero boundary value data on the interior of
$C_1\cup C_2.$
 Actually,
let $C_3$ be an arc of geodesic lying in $\Delta$ joining a point
$c_1$ on  $C_1$ to a point $c_2$ on $C_2.$ Let $a=C_1\cap C_2$ and
let $\Delta_0$ be the geodesic triangle  with vertices $a, c_1$
and $c_2.$

Let   $f$ be the restriction of $\phi$ to $C_3.$ Notice that the
solution of the Dirichlet problem on $\Delta_0$ taking zero
boundary value data on the sides $[a, c_1], [a, c_2]$ and   $f$ on
the side $[c_1, c_2]$ gives rise to a superior barrier at any
point of the interior of $C_1\cup C_2.$ Of course the zero
function is an inferior barrier to the problem. Thus $u$ takes the
desired boundary value data, as we claimed.}
\end{example}
%Uniqueness follows from the general maximum principle \cite{HRS}.

\begin{example}\label{scherkideal}{\em

We consider a geodesic triangle  $\Delta$ in $\hip^2\times\{0\}$
with two vertices $a, b$ on $\pain\hip^2\times\{0\}$  and a third
vertex $c$ on $\hip^2\times\{0\}$. 
 Doing a similar construction
as in \exaref{scherk}, we can solve our Dirichlet problem on
$\Delta$ taking infinite boundary value data on the complete
geodesic $(a, b)$ and zero boundary value data on the two other
sides.

Assume now that the interior angle at vertex $c$
is $\pi/k$, $k\in \Na^\ast$. Using  Schwarz reflection on the geodesics arcs
$(c, a),$
 $(c, b),$ and successively about the geodesic boundaries as well, we
obtain a complete embedded minimal surface  in $\hip^2\times \R$.
Since the angle at the vertex $c$ is
$\pi/k,\, k\in\Na^\ast,$ we get a complete graph over  an ideal
geodesic polygon with $2k$ sides,  taking successively boundary
values $+\infty$ and $-\infty$. These  minimal complete
graphs  can also be built combining some results on harmonic maps
from the complex plane into the hyperbolic plane, done in \cite
{T-W}, \cite{H-T-T-W}  and \cite{hst}.
We observe that these examples are a particular case of a general
result found in \cite{C-R}.}
\end{example}

\section{Minimal graphs with finite and asymptotic boundary in $\hd \times \R$}
\label{mai}

\begin{lemma} \label{l1}
 Let $\rho >0$ and let $\Cc_\rho\subset \hd\times \{0\}$ be a circle
 of radius $\rho.$ Then there exists a unique catenoid $\cro$ in
 $\hip^2\times\R$
   orthogonal to the slice $\hd\times
 \{0\}$ along $\Cc_\rho.$
 Its asymptotic boundary is $\pain \hd\times \{\pm t_0\},$ for some
 $0<t_0,$ where   $t_0:=f(\rho)$ is  an increasing function of $\rho$  given by

\begin{multline}\label{f1}
 \hfill f(\rho)=\int\limits_{\rho}^\infty\frac{\sinh\rho}{\sqrt{\sinh^2
r-\sinh^2 \rho}}\, \rmd r\hfill
\end{multline}

Furthermore,  $\lim\limits_{\rho\goto 0} f(\rho)=0, $ and 
$\lim\limits_{\rho\goto \infty} f(\rho)=\pi/2. $
\end{lemma}

The proof of  \lemref{l1} follows from  Proposition 5.1 of
\cite{N-SE-S-T} and \cite{SE-T1}. For later use we call $\cro^+$
(resp. $\cro^-$) the part of the catenoid $\cro$ in $\hip^2\times
[0, \infty)$(resp.  in $\hip^2\times (-\infty, 0]$).

\begin{proposition}[A characterization of minimal vertical
graphs]\label{p1}
 Let $M$ be a minimal surface immersed in
$\hip^2\times \R,$ whose finite boundary is a Jordan curve $\Ga$ and
whose asymptotic boundary is $\pain \hd\times \{t_0\}\subset \pain
\hd \times \R,\; t_0\geqs 0.$ Assume that $\Ga$ is a vertical
graph over a Jordan curve $C\subset \hip^2\times \{0\}.$ Assume
also that the vertical
projection of $M$ is contained in $\ext C.$\\
 Then $M$ is a vertical graph. Furthermore, if $\Ga=C, $ then $t_0<\pi/2$ and $M$
 inherits all symmetries of $\Ga.$ Particularly, if $\Ga$ is an
 horizontal circle then $M$ is part of a catenoid.

\end{proposition}

\begin{proof} The proof is somewhat straightforward. We will
just sketch it as follows. The first statement is a consequence of
Alexandrov Reflection Principle on  horizontal  slices doing
vertical reflections. The second statement ($t_0<\pi/2$)  is a
consequence of \lemref{l1} using the family of catenoids, coming
from the infinity towards $M$. The third statement is a
consequence of Alexandrov reflection Principle on vertical
geodesic planes.

\end{proof}

The following Remark is inferred from \lemref{l1} and maximum
principle.

\begin{remark}\label{r2}{\em  Let $C_\rho$ be a circle of radius
$\rho$ in $\hip^2\times \{0\},$  and let  $\Gai\subset \pain \hd
\times \R,$ be a Jordan curve that is a vertical graph over $\pain
\hip^2\times\{0\}.$ If the height function $t$ of $\Gai$ satisfies $ t
>f(\rho)$ there is no minimal vertical graph over $\ext(C_\rho)$
whose finite boundary is $C_\rho$ and whose asymptotic boundary is
$\Gai.$}
\end{remark}

\begin{definition}[Admissible unbounded domains in
$\hip^2$]\label{d2}{\em 
 Let $\Om$ be an unbounded domain in the slice
$\hip^2\times\{0\}$ and let $\partial \Om$ be its boundary.  We
say that $\Om$ is an {\em admissible domain} if each connected component
$C_0$ of $\partial \Om$ satisfies one of the following conditions:
\begin{enumerate}
\item $C_0$ is a Jordan curve.
\item $C_0$ is a properly embedded curve  such that the asymptotic
boundary is one point.
\item $C_0$ is a properly embedded curve  such that the asymptotic
boundary is two distinct points.
\end{enumerate}

Finally, each connected component $C_0$ of $\partial \Om$
satisfies the {\em Exterior circle of $($uniform$)$ radius $\rho$
condition}, that is,  at any point $p\in C_0$ there exists a
circle $C_\rho$ of radius $\rho$ such that $p\in C_0\cap C_\rho$
and $\overl{\inter
 C_\rho}\cap \Om=\emptyset.$

 If $\Om$ is an unbounded admissible domain then  we denote
 by $\rho_\Om$   the supremum of the set of these $\rho.$

If the components of $\partial\Om$ are compact, we set
 $C:=\partial\Om,$ hence  $C=C_1\cup\ldots\cup C_n$ is the union of disjoint
Jordan curves $C_j,\, j=1,\ldots,n $ with pairwise disjoint
interiors.   We
 set $\ext (C)=\ext(C_1)\cap\cdots \cap\ext(C_n)$ and
 $\inter(C)=\inter(C_1)\cup\cdots \cup\inter(C_n).$ In this case
 we set $\rho_C:=\rho_\Om.$}

\end{definition}

In the next theorem we need the function  $f(\rho)$   given in
\lemref{l1} ({\em height} of the {\em catenoid} $\Mc_\rho$ arising
orthogonally from the slice along a circle of radius $\rho_\Om$).

\begin{theorem}\label{t1}
Let $\Om$ be an admissible unbounded domain.  Let \\
$g: \partial\Om\cup\partial_\infty\Om \rightarrow \R$ be a  
bounded function 
 taking zero  boundary value data on $\partial\Om$,
everywhere continuous except maybe at a finite set
$S\subset \partial_\infty\Om$. Let
$\Gai\subset \pain \hd \times \R$ be the union of the graph of $g$
restricted to $\partial_\infty\Om $ with the vertical segments at
the points of $\partial_\infty\hip^2$ of discontinuities of $g.$
 \\
If   the height function  $t$ of $\Gai$ satisfies  $-f(\rho_\Om)\leqs
t\leqs f(\rho_\Om),$ then there exists a minimal vertical graph
over $\Om$ with finite boundary $\partial \Om$ and asymptotic
boundary $\Gai.$
%Particularly, if $\Gai$ is a circle we obtain a minimal graph with
%boundary $C$ asymptotic to a catenoid.

\vskip1mm
 Particularly, if $C\subset\hip^2\times \{0\}$ is  
a Jordan curve satisfying the 
{\em Exterior circle of radius $\rho$ condition} and if    $g:
\partial_\infty\hip^2 \rightarrow \R$ is  a continuous function
satisfying 
$-f(\rho_\Om)\leqs g(p) \leqs f(\rho_\Om) $   at any point $p\in
\pain \hd,$  then there exists a unique minimal vertical graph
over $\Om$ with finite boundary $\partial \Om$ and asymptotic
boundary $\Gai$.

Finally, there is no such minimal
 graph, if  $\partial \Om$ is compact and the height function $t$  of $\Gai $
satisfies $\vert t \vert >\pi/2$.

%Particularly, if $\Gai$ is a circle we obtain a minimal graph with
%boundary $C$ asymptotic to a catenoid.

\end{theorem}

\begin{proof}
 Consider the family of catenoids $\Mc_\rho$ given by \lemref{l1}.
Notice that our assumptions imply that at each point $p\in C$
there exists a circle $\Cc_{\rho_C}(p)$ of radius $\rho_C$
contained in $\overl{\inter(C)}$ with $p\in C_\rho(p)\cap C.$ Let
$\Mc_{\rho_C}^+(p)$ and $\Mc_{\rho_C}^-(p)$   be the the upper and
lower half- catenoids cutting orthogonally the slice  $t=0$ along
the circle $\Cc_{\rho_C}(p)$.

Take one of these lower half-catenoids as a subsolution, and take
one of these upper half-catenoids as the supersolution $\phi$ in
Perron process, \thmref{perron}. It follows that this family of
half- catenoids provide also a family of  barriers at each point
of $C$ to our problem in the sense of \defref{d1}-(\ref{D.Item.bar.2}).
Therefore
our Dirichlet Problem (P), see \defref{probP}, can be solved 
using \coref{C.bounded.data}.

  If $C\subset\hip^2\times \{0\}$ is a Jordan curve
satisfying the {\em Exterior circle of radius
$\rho$ condition},  and if   $g:
\partial_\infty\hip^2 \rightarrow \R$ is   continuous, then the uniqueness follows from the classical
maximum principle \thmref{mp}. This proves the first assertion of
the statement.

\vskip2mm
 To prove the nonexistence part
assume by contradiction that there exists a solution $u$ such that
the height function $t$ of  $\Gai $ satisfies $t >\pi/2$. Notice  that the
graph of $u$ is above the slice $t=0.$ Now choose a catenoid
$\Mc_\rho$ with $t$ axis and large ``neck'' ($\rho$ big enough)
disjoint from the graph of $u.$ Let $\Mc_\rho(\ep)= \Mc_\rho +
\ep$  be the $\ep$-vertical translation of $\Mc_\rho$, with
$\ep>0$ small enough. Now shrink the catenoid $\Mc_\rho(\ep)$ in
the family of catenoids with the same axis making the ``neck''
going to zero. We will find a first interior point of 
contact of the graph of 
$u$ with one of these catenoids. This
gives a contradiction by the maximum principle and completes the
proof of the Theorem.

\end{proof}

\begin{remark}{\em
 A computation shows that any catenoid in $\hd \times \R$ has finite total
extrinsic curvature. We set here a question: is it true that the same holds for
any
exterior minimal graph in $\hd \times \R$?}

\end{remark}

We now  restrict our attention to certain admissible
domains such that each component of the boundary has two points at
its asymptotic boundary and has at each point of its finite
boundary an exterior equidistant curve. To be more precise:

\begin{definition}[E-admissible unbounded domains in
$\hip^2$]\label{d3}{\em
 Let $\Om$ be an unbounded domain in the slice
$\hip^2\times\{0\}$ and let $\partial \Om$ be its boundary.  We
say that $\Om$ is an {\em E-admissible domain} if

\begin{enumerate}

\item \label{D.admiss.1} Each connected component $C_0$ of $\partial \Om$ is a
properly
embedded curve such that the asymptotic boundary is constituted of
two distinct points.

 \item  We require that there exists $r>0$ such that  each point of $\partial \Om$
satisfies the {\em Exterior equidistant curve of $($uniform$)$
curvature $\tanh r$ condition}; that is,  at any point $p\in
\partial \Om$ there exists an equidistant curve  $E_{r}$ of
curvature $\tanh r$ (with respect to the exterior unit normal to
$\Om$ at $p$ ), with $p\in
\partial\Om\cap E_r$ and $ E_r\cap \Om=\emptyset.$

\end{enumerate}
}
\end{definition}
Thus every E-admissible domain is an admissible domain.

If $\Om$ is a convex domain satisfying the condition (\ref{D.admiss.1})
of \defref{d3} then $\Om$ is an E-admissible domain.

If each connected component $C_0$ of $\partial \Om$ is  an
equidistant curve then $\Om$ is an E-admissible (maybe nonconvex)
domain.

If $\Om$ is an unbounded E-admissible domain then  we denote
 by $r_\Om$ $\geqs 0$   the infimum  of the set of these $r.$
 If $\Om$ is a convex E-admissible domain then $r_\Om=0$.

We will use in the next result the function $H$ defined by 
Formula \eqref{Eq.function.H}  in 
Proposition \ref{P.hyperbolic translation}.

\begin{theorem}\label{t2}

Let $\Om$ be an  E-admissible unbounded domain.  Let \\
$g: \partial\Om\cup\partial_\infty\Om \rightarrow \R$ be a  
bounded function 
 taking zero  boundary value data on $\partial\Om$,
everywhere continuous except maybe at a finite set
$S\subset \partial_\infty\Om$.
Let
$\Gai\subset \pain \hd \times \R$ be the union of the graph of $g$
restricted to $ \partial_\infty\Om$ with the vertical segments at
the points of discontinuities of $g.$
 \\
If   the height function  $t$ of $\Gai$ satisfies  
$-H(\cosh r_\Om)\leqs t \leqs H(\cosh r_\Om)$,
 then there exists a minimal vertical graph over $\Om$
with finite boundary $\partial \Om$ and asymptotic boundary
$\Gai.$

\end{theorem}

\begin{proof}
The proof is the same as in \thmref{t1}, replacing the minimal
catenoids by the  minimal surfaces invariant by hyperbolic
translations $M_d,\ d>1$, given in Proposition \ref{P.hyperbolic translation}.
This completes the proof of the Theorem.
\end{proof}

Notice that if $\Om$ is convex then it is E-admissible and
$r_\Om=0,$ thus $H(\cosh r_\Om)=\infty.$

\end{document}